\documentclass[11pt, reqno]{amsart}

\usepackage{amsthm,amssymb,amstext,amscd,amsfonts,amsbsy,amsrefs,amsxtra,latexsym,amsmath,xcolor,mathrsfs,fancybox,upgreek, soul,url}
\usepackage[english]{babel}
\usepackage[all,cmtip]{xy}
\usepackage[latin1]{inputenc}
\usepackage{cancel}
\usepackage[draft]{hyperref}
\usepackage{comment,enumitem}
\usepackage{mdframed}
\allowdisplaybreaks
\usepackage{mathtools}
\usepackage{enumerate}

\newtheorem{theorem}{Theorem}

\theoremstyle{definition}

\theoremstyle{remark}

\newtheorem*{remark}{Remark}

\numberwithin{theorem}{section}
\numberwithin{equation}{section}

\newcommand{\N}{\mathbb{N}}
\newcommand{\Z}{\mathbb{Z}}
\newcommand{\R}{\mathbb{R}}
\newcommand{\C}{\mathbb{C}}

\newcommand{\Q}{\mathbb{Q}}

\newcommand{\sgn}{\operatorname{sgn}}

\newcommand{\J}{\mathscr{S}}
\def\H{\mathbb{H}}

\begin{document}

\title[A short note on higher Mordell integrals]{A short note on higher Mordell integrals}

\author{Joshua Males}

\address{Department of Mathematics and Computer Science, Division of Mathematics, University of Cologne, Weyertal 86-90, 50931 Cologne, Germany}
\email{jmales@math.uni-koeln.de}

\begin{abstract}
	Classical mock modular and quantum modular forms are known to have an intimate relationship with Mordell integrals thanks to Zwegers groundbreaking Ph.D.\@ thesis. More recently, generalisations of mock/quantum modular forms to so-called ``higher depth" versions have been intensively studied. In essence, a mock/quantum modular form of depth $d$ is such that the error of modularity transforms as another mock/quantum modular form of depth $d-1$. In this short note we use techniques of Bringmann, Kaszian, and Milas to show that the double Eichler integrals of a family of depth two quantum modular forms of weight one previously studied by the author can be related to certain ``higher" Mordell integrals, meaning it may be written as a certain double integral, \`{a} la Zwegers.
\end{abstract}

\maketitle

\section{Introduction}
The Mordell integral
\begin{equation}\label{Equation: one-dim Mordell}
h(z) = h(z;\tau) \coloneqq \int\limits_{\R} \frac{\cosh(2 \pi z w)}{\cosh(\pi w)} e^{\pi i \tau w^2} dw,
\end{equation}
where $z \in \C$ and $\tau \in \H$, is intricately linked to various areas of number theory. In particular, classical results show the connection between specialisations of \eqref{Equation: one-dim Mordell} are connected to the Riemann zeta function \cite{siegel1932uber}, Gauss sums \cite{kronecker1889bemerkungen,kronecker1889summirung}, and class number formulas \cite{mordell1920,mordell1933}. 

More recently, Zwegers used Mordell integrals to describe the completion of Lerch sums in his celebrated thesis \cite{zwegers2008mock}. In particular, Zwegers observed that we can relate \eqref{Equation: one-dim Mordell} to an Eichler integral in the following way
\begin{equation}\label{Equation: h(a tau - b)}
h(a\tau - b) = -e^{-2 \pi i a\left(b+ \frac{1}{2}\right)} q^{\frac{a^2}{2}} \int_{0}^{i \infty} \frac{g_{a+ \frac{1}{2} , b+ \frac{1}{2}} (w)}{\sqrt{-i (w + \tau)}} dw.
\end{equation}
Here, $g_{a,b}$ is the weight $\frac{3}{2}$ unary theta function given by ($a,b \in \R$)
\begin{equation*}
g_{a,b} (\tau) \coloneqq \sum_{n \in a+\Z} ne^{2 \pi i bn} q^{\frac{n^2}{2}}.
\end{equation*}
Zwegers then showed that a modular completion of Lerch sums may be found. To do so, he found that the error of modularity $h(a\tau -b)$ also appears when considering integrals of the same form as \eqref{Equation: h(a tau - b)} with lower integration boundary  $-\bar{\tau}$ instead of $0$.

Furthermore, Eichler integrals of the form
\begin{equation*}
\int_{-\bar{\tau}}^{i \infty} \frac{g (w)}{(-i (w + \tau))^{\frac{3}{2}}} dw
\end{equation*}
with $g$ a cuspidal theta function have been studied by many authors in recent times, perhaps most notably in relation to quantum modular forms, e.g. \cite{bringmann2016half,folsom2017strange,rolen2013strange}. Quantum modular forms were introduced by Zagier in \cite{zagier2001vassiliev,zagier2010quantum} and are essentially functions $f \colon \mathcal{Q} \rightarrow \C$ for some fixed $\mathcal{Q} \subseteq \mathbb{Q}$, whose errors of modularity (for $M = \left( \begin{smallmatrix} a & b \\ c & d \end{smallmatrix} \right) \in \Gamma \subset \text{SL}_2 (\Z)$)
\begin{equation*}
f(\tau) - (c \tau + d)^k f (M \tau)
\end{equation*}
are in some sense ``nicer" than the original function. Often, for example, the original function $f$ is defined only on $\Q$, but the ``errors of modularity" can be defined on some open subset of $\R$. Quantum modular forms have been the topic of much interest in the past decade,  for example there is a fascinating connection between them and mock modular forms - surveyed in \cite{ono2009unearthing} - which has been investigated in papers such as \cite{bringmann2015unimodal,bringmann2016half,bryson2012unimodal}, among others. Interesting examples of quantum modular forms also lie at the interface of physics and knot theory, see e.g. a study of Kashaev invariants of $(p,q)$-torus knots in \cite{hikami2003torus,hikami2015torus} and investigations of Zagier into limits of quantum invariants of $3$-manifolds and knots \cite{zagier2010quantum}. Understanding the error of modularity of quantum modular forms is then clearly an important problem. This paper serves to extend results of Bringmann, Kaszian, and Milas to a certain infinite family of so-called quantum modular forms of depth two (see the sequel for precise definitions).

 A certain generalisation of quantum modular forms was introduced in \cite{HigherDepthQMFs}. The authors define so-called higher depth quantum modular forms, and provide two examples of such forms of depth two that arise from characters of vertex operator algebras. In the simplest case, quantum modular forms of depth two are functions that satisfy
 \begin{equation*}
 f(\tau) - (c \tau + d)^k f(M \tau) \in \mathsf{Q}_k (\Gamma) \mathcal{O} (R) + \mathcal{O} (R),
 \end{equation*}
 where $\mathsf{Q}_k (\Gamma)$ is the space of quantum modular forms of weight $k$ on $\Gamma$, and $\mathcal{O}(R)$ is the space of real-analytic functions on $R \subset \R$. A crucial step in showing the generalised quantum modularity property of their functions $F_1$ and $F_2$ is the appearance of a two-dimensional Eichler integral of the shape
\begin{equation}\label{Equation: form of double Eichler integral}
\int_{- \bar{\tau}}^{i \infty} \int_{\omega_1}^{i \infty} \frac{ g_1(\omega_1) g_2(\omega_2) }{ \sqrt{-i(\omega_1 + \tau)} \sqrt{-i(\omega_2 + \tau)}} d \omega_2 d \omega_1 ,
\end{equation}
where the $g_j $ lie in the space of vector-valued modular forms on $\text{SL}_2 (\Z)$. In the next paper in the series of Bringmann, Kaszian, and Milas \cite{HigherDepthQMFs2} the connection between such a two-dimensional Eichler integral and higher Mordell integrals is explored, in particular with the example of the function $F_1$ carried over from \cite{HigherDepthQMFs}. 

In particular, the higher Mordell integrals investigated in \cite{HigherDepthQMFs2} provide the error of modularity of their function $F_1$, in turn developing the theory of Zwegers to higher dimensions. In the present paper, we take the example of depth two Mordell integrals of \cite{HigherDepthQMFs2} and extend this to an infinite family of similar functions, thereby also providing an infinite family of errors of modularity of the relatively new higher depth quantum modular forms. One therefore also sees that by understanding higher Mordell integrals, we already obtain intrinsic information about the higher depth quantum modular form. Furthermore, the construction given in this paper gives hints as to how one could obtain similar results for arbitrary depth. The outline is sketched in the following.

In \cite{2018arXiv181001341M} a family of functions is given as a generalisation of the function $F_1$. Each function $F$ in this more general family from \cite{2018arXiv181001341M} is of the shape (up to addition by one-dimensional partial theta functions) 
\begin{equation*}
 \sum_{\alpha \in \J} \varepsilon(\alpha) \sum_{n \in \alpha + \N_0^2 } q^{Q(n)},
\end{equation*} 
with $Q(x)= a_1 x_1^2 + a_2 x_1 x_2 + a_3 x_2^2$ a positive definite integral binary quadratic form, $\J$ a finite set of pairs $\alpha \in \Q^2 \backslash\{(0,0)\}$, and $\varepsilon \colon \J \rightarrow \R \backslash \{0\}$. Each $F$ is shown to be vector-valued quantum modular form of depth two and weight one. Similarly to \cite{HigherDepthQMFs}, a key compenent is the introduction of the double Eichler integral
\begin{equation}\label{Equation: definition of E (tau)}
\mathcal{E}_{\alpha} (\tau) \coloneqq - \frac{\sqrt{D}}{4} \int_{- \bar{\tau}}^{i \infty} \int_{\omega_1}^{i \infty} \frac{\theta_1 (\alpha; \omega_1, \omega_2) + \theta_2 (\alpha; \omega_1, \omega_2)}{\sqrt{-i (\omega_1 + \tau)} \sqrt{-i(\omega_2 + \tau)}} d\omega_2 d\omega_1,
\end{equation}
where $D \coloneqq 4a_1a_3 - a_2^2 > 0$, and $\theta_1, \theta_2$ are given explicitly in Section \ref{Section: Mordell integrals}. It is shown in \cite{2018arXiv181001341M} that using Shimura theta functions we may rewrite this in the form \eqref{Equation: form of double Eichler integral}.

This short note serves to show that techniques of Bringmann, Kaszian, and Milas of relating their double Eichler integral to higher Mordell integrals in \cite{HigherDepthQMFs2} immediately carry over to the more general setting of \cite{2018arXiv181001341M}. In a similar fashion to \cite{HigherDepthQMFs2} we define
\begin{equation*}
H_\alpha (\tau) \coloneqq - \sqrt{D} \int_{0}^{\infty} \int_{\omega_1}^\infty \frac{ \theta_1 (\alpha; \omega_1, \omega_2) + \theta_2 (\alpha; \omega_1, \omega_2)}{  \sqrt{-i (\omega_1 + \tau)} \sqrt{-i (\omega_1 + \tau)} } d\omega_1 d\omega_2,
\end{equation*}
along with the functions

\begin{equation*}
\mathcal{F}_\alpha (x) \coloneqq \frac{\sinh(2 \pi x)}{\cosh(2 \pi x) - \cos(2 \pi \alpha)}, \hspace*{20pt}     \mathcal{G}_\alpha(x)\coloneqq \frac{\sin(2 \pi \alpha)}{\cosh(2 \pi x) - \cos(2 \pi \alpha)}.
\end{equation*}
Our result is the following theorem (there is also a related expression for $\alpha \in \Z^2$ , taking first a limit in $\alpha_1$ and using the same method as below, and then taking a limit in $\alpha_2$).

\begin{theorem}\label{Theorem: main}
	For $\alpha \not\in \Z^2$, we have that 
	
	\begin{equation*}
	H_\alpha (\tau) = \int_{\R^2} g_\alpha(\tau) d\omega,
	\end{equation*}
	where we set
	
	\begin{equation*}
	g_\alpha(\tau) \coloneqq
	\begin{dcases}
	2 \mathcal{G}_{\alpha_1}(\omega_1) \mathcal{G}_{\alpha_2} (\omega_2) - 2 \mathcal{F}_{\alpha_1} (\omega_1) \mathcal{F}_{\alpha_2} (\omega_2) & \text{ if } \alpha_1, \alpha_2 \not\in \Z, \\
	 - 2 \mathcal{F}_{0} (\omega_1) \mathcal{F}_{\alpha_2} (\omega_2) + \frac{2}{\pi \omega_1} \mathcal{F}_{\alpha_2} \left( \omega_2 + \frac{a_2}{2a_3} \omega_1  \right) & \text{ if } \alpha_1 \in \Z, \alpha_2 \not\in \Z, \\
	 - 2  \mathcal{F}_{\alpha_1} (\omega_1) \mathcal{F}_{0} (\omega_2) + \frac{2}{\pi \omega_2} \mathcal{F}_{\alpha_1} \left( \omega_1 + \frac{a_2}{2a_1} \omega_2  \right) & \text{ if } \alpha_1 \not\in \Z, \alpha_2 \in \Z. \\
	\end{dcases}
	\end{equation*}
\end{theorem}

\section{Preliminaries}
Here we recall a few relevant results on double error functions that we need in the rest of this note. We first define a rescaled version of the usual one-dimensional error function. For $u \in \R$ set
\begin{equation}
E(u) \coloneqq 2 \int_0^u e^{- \pi \omega^2} d \omega .
\end{equation}
We also require, for non-zero $u$, the function
\begin{equation*}
M(u) \coloneqq \frac{i}{\pi} \int_{\R - iu} \frac{e^{ - \pi \omega^2 - 2 \pi i u \omega}}{\omega} d \omega.
\end{equation*}
A relation between $M(u)$ and $E(u)$, for non-zero $u$, is given by
\begin{equation}
M(u) = E(u) - \sgn(u).
\end{equation}
We further need the two-dimensional analogues of the above functions. Following \cite{GeneralisedErrorFunctions} and changing notation slightly, we define $E_2 \colon \R \times \R^2 \rightarrow \R$ by
\begin{equation*}
E_2 (\kappa ; u) \coloneqq \int_{\R^2} \sgn(\omega_1) \sgn(\omega_2 + \kappa \omega_1) e^{- \pi \left( (\omega_1 - u_1)^2 + (\omega_2 - u_2)^2 \right)} d \omega_1 d \omega_2, 
\end{equation*}
where throughout we denote components of vectors just with subscripts. Again following \cite{GeneralisedErrorFunctions}, for $u_2, u_1 - \kappa u_2 \neq 0$, we define
\begin{equation}\label{Equation: definition of M2}
M_2 (\kappa; u_1, u_2) \coloneqq - \frac{1}{\pi^2} \int_{\R - iu_2} \int_{\R - i u_1} \frac{ e^{ - \pi \omega_1^2 - \pi \omega_2^2 - 2 \pi i (u_1 \omega_1 + u_2 \omega_2)}}{\omega_2 (\omega_1 - \kappa \omega_2)} d\omega_1 d \omega_2 .
\end{equation}
Then we have that
\begin{equation}\label{Equation: relation between M_2 and E_2}
\begin{split}
M_2 (\kappa; u_1 , u_2) = & E_2(\kappa; u_1 , u_2) - \sgn(u_2) M (u_1) \\
& - \sgn(u_1 - \kappa u_2) M \left( \frac{u_2 + \kappa u_1}{\sqrt{1 + \kappa^2}} \right) - \sgn(u_1) \sgn(u_2 + \kappa u_1).
\end{split}
\end{equation}
The relation \eqref{Equation: relation between M_2 and E_2} extends the definition of $M_2 (u)$ to include $u_2 = 0$ or $u_1 = \kappa u_2$ - note however that $M_2$ is discontinuous across these loci. Further, it is shown in the proof of Lemma 7.1 of \cite{2018arXiv181001341M} that for $u = u(n) \coloneqq (2 \sqrt{a_1} n_1 + \frac{a_2}{\sqrt{a_1}} n_2,  m n_2)$, along with $\kappa \coloneqq \frac{a_2}{\sqrt{D}}$, and $
 m \coloneqq \sqrt{4a_3 - \frac{a_2^2}{a_1}}$ we have that
\begin{equation}\label{Equation: integral form of M2}
\begin{split}
M_2 (\kappa; \sqrt{v} u)  = & - \frac{\sqrt{D} n_2 (2 a_1 n_1 + a_2 n_2)}{2 a_1} q^{Q(n)} \int_{- \bar{\tau}}^{i \infty} \frac{ e^{\frac{\pi i (2 a_1 n_1 + a_2 n_2)^2 \omega_1}{2 a_1}}}{\sqrt{-i (\omega_1 + \tau)}} \int_{\omega_1}^{i \infty} \frac{e^{\frac{\pi i D n_2^2 \omega_2}{2 a_1}}}{\sqrt{-i (\omega_2 + \tau)}} d\omega_2 d\omega_1\\
& - \frac{\sqrt{D} n_1 (a_2 n_1 + 2 a_3 n_2)}{2 a_3} q^{Q(n)}  \int_{- \bar{\tau}}^{i \infty} \frac{e^{\frac{\pi i (a_2 n_1 + 2 a_3 n_2)^2 \omega_1}{2 a_3}}}{\sqrt{-i (\omega_1 + \tau)}} \int_{\omega_1}^{i \infty} \frac{e^{\frac{\pi i D n_1^2 \omega_2}{2 a_3}}}{\sqrt{-i (\omega_2 + \tau)}} d\omega_2 d\omega_1.
\end{split}
\end{equation}

\section{Proof of Theorem \ref{Theorem: main}}\label{Section: Mordell integrals}
\begin{proof}
	By analytic continuation it suffices to show that the theorem holds for $\tau = iv$, and we begin by showing that
\begin{equation*}
H_\alpha (iv) = 2 \lim\limits_{r \rightarrow \infty} \sum_{\substack{n \in \alpha + \Z^2 \\ |n_j - \alpha_j | \leq r}} M_2 \left(\kappa; \sqrt{\frac{v}{2}} u \right) e^{2 \pi v Q(n)},
\end{equation*}
We begin with the expression \eqref{Equation: integral form of M2} evaluated at $\tau = iv$, giving
\begin{equation}\label{Equation: M2 evaluated at tau = iv}
\begin{split}
M_2 (\kappa; \sqrt{v} u)  = & - \frac{\sqrt{D} n_2 (2 a_1 n_1 + a_2 n_2)}{2 a_1} q^{Q(n)} \int_{iv}^{i \infty} \frac{ e^{\frac{\pi i (2 a_1 n_1 + a_2 n_2)^2 \omega_1}{2 a_1}}}{\sqrt{-i (\omega_1 + iv)}} \int_{\omega_1}^{i \infty} \frac{e^{\frac{\pi i D n_2^2 \omega_2}{2 a_1}}}{\sqrt{-i (\omega_2 + iv)}} d\omega_2 d\omega_1\\
& - \frac{\sqrt{D} n_1 (a_2 n_1 + 2 a_3 n_2)}{2 a_3} q^{Q(n)}  \int_{iv}^{i \infty} \frac{e^{\frac{\pi i (a_2 n_1 + 2 a_3 n_2)^2 \omega_1}{2 a_3}}}{\sqrt{-i (\omega_1 + iv)}} \int_{\omega_1}^{i \infty} \frac{e^{\frac{\pi i D n_1^2 \omega_2}{2 a_3}}}{\sqrt{-i (\omega_2 + iv)}} d\omega_2 d\omega_1.
\end{split}
\end{equation}
We make the shift $\omega_j \rightarrow 2 i \omega_j + iv$. The terms in the exponential in the first term on the right-hand side become 
\begin{equation*}
\frac{\pi i (2 a_1 n_1 + a_2 n_2)^2 (2 i \omega_1 + iv)}{2 a_1} = - \pi \frac{(2 a_1 n_1 + a_2 n_2)^2}{a_1}  \omega_1 - \pi v \frac{ (2 a_1 n_1 + a_2 n_2)^2 }{2 a_1},
\end{equation*}
along with
\begin{equation*}
\frac{\pi i D n_2^2 (2 i \omega_2 + iv)}{2 a_1} = - \pi \frac{ D n_2^2}{a_1}  \omega_2 - \pi v \frac{ D n_2^2 }{2 a_1}.
\end{equation*}
Pulling out the above two terms dependent on $v$ gives $- 2 \pi v Q(n)$. Then we see that the first term on the right-hand side of \eqref{Equation: M2 evaluated at tau = iv} is equal to
\begin{equation*}
\begin{split}
 \frac{ \sqrt{D} n_2 (2 a_1 n_1 + a_2 n_2)}{a_1} e^{- 4 \pi v Q(n)} \int_{0}^{\infty} \frac{ e^{-\frac{\pi (2 a_1 n_1 + a_2 n_2)^2 \omega_1}{ a_1}}}{\sqrt{\omega_1 + v}} \int_{\omega_1}^{\infty} \frac{e^{- \frac{\pi D n_2^2 \omega_2}{a_1}}}{\sqrt{\omega_2 + v}} d\omega_2 d\omega_1.
\end{split}
\end{equation*}
A similar expression holds for the second term, and thus we can write $e^{4 \pi v Q(n)} M_2 (\kappa; \sqrt{v}u) $ as the sum of the two terms 
\begin{equation*}
\begin{split}
 \frac{ \sqrt{D} n_2 (2 a_1 n_1 + a_2 n_2)}{a_1} \int_{0}^{\infty} \int_{\omega_1}^{\infty} \frac{ e^{- \frac{\pi (2 a_1 n_1 + a_2 n_2)^2 \omega_1}{a_1} - \frac{\pi D n_2^2 \omega_2}{a_1}}}{\sqrt{\omega_2 + v} \sqrt{\omega_1 + v}} d\omega_2 d\omega_1,
\end{split}
\end{equation*}
and
\begin{equation*}
\begin{split}
\frac{ \sqrt{D} n_1 (a_2 n_1 + 2 a_3 n_2)}{a_3} \int_{0}^{\infty} \int_{\omega_1}^{\infty} \frac{ e^{-\frac{\pi (a_2 n_1 + 2 a_3 n_2)^2 \omega_1}{a_3} - \frac{\pi D n_1^2 \omega_2}{a_3}}}{\sqrt{\omega_2 + v} \sqrt{\omega_1 + v}} d\omega_2 d\omega_1.
\end{split}
\end{equation*}
Let $v \rightarrow \frac{v}{2}$, sum over $n \in \alpha + \Z^2$ such that $|n_j - \alpha_j | \leq r$ and let $r \rightarrow \infty$. In the same way as \cite{HigherDepthQMFs2} we may use Lebesgue's dominated convergence theorem to obtain
\begin{equation}\label{Equation: final step in showing link between double error and Mordell integral}
2 \lim\limits_{r \rightarrow \infty} \sum_{\substack{n \in \alpha + \Z^2 \\ |n_j - \alpha_j | \leq r}} M_2 \left(\kappa; \sqrt{\frac{v}{2}} u \right) e^{2 \pi v Q(n)} = - \sqrt{D} \int_{0}^{\infty} \int_{\omega_1}^{\infty}  \frac{\theta_1 (\alpha; \omega) + \theta_2(\alpha; \omega)}{{\sqrt{\omega_2 + iv} \sqrt{\omega_1 + iv}}} d\omega_2 d\omega_1,
\end{equation}
where we set
\begin{equation*}
\theta_1 (\alpha; \omega_1, \omega_2) \coloneqq \frac{1}{a_1} \sum_{n \in \alpha + \Z^2} (2 a_1 n_1 + a_2 n_2) n_2 e^{\frac{\pi i (2 a_1 n_1 + a_2 n_2)^2 \omega_1}{2 a_1} + \frac{ \pi i D  n_2^2 \omega_2}{2 a_1}}
\end{equation*}
and
\begin{equation*}
\theta_2 (\alpha; \omega_1, \omega_2) \coloneqq \frac{1}{a_3} \sum_{n \in \alpha + \Z^2} ( a_2 n_1 + 2 a_3 n_2) n_1 e^{\frac{\pi i (a_2 n_1 + 2 a_3 n_2)^2 \omega_1}{2 a_3} +\frac{ \pi i D n_1^2 \omega_2}{2 a_3}}.
\end{equation*}
Further, it is clear by definition that the right-hand side of \eqref{Equation: final step in showing link between double error and Mordell integral} is equal to $H_\alpha(iv)$, and so we have shown the first claim.

\begin{remark}
	We note that these theta functions are exactly those appearing in the double Eichler integral associated to the family of quantum modular forms of depth two given in \cite{2018arXiv181001341M}.
\end{remark}

Now we concentrate on $M_2(\kappa; \sqrt{v}u)$. Assuming that $u_2 , u_1 - \kappa u_2 \neq 0$ (which happens precisely when $\alpha_1,\alpha_2 \not\in \Z$), rewriting \eqref{Equation: definition of M2} implies that
\begin{equation*}
M_2(\kappa; \sqrt{v}u) = - \frac{1}{\pi^2} e^{- \pi v (u_1^2 + u_2^2)} \int_{\R^2} \frac{e^{-\pi v \omega_1^2 - \pi v \omega_2^2}}{(\omega_2 - iu_2) (\omega_1 - \kappa \omega_2 - i(u_1 - \kappa u_2))} d\omega_1 d\omega_2.
\end{equation*}
Plugging in our definition of $u(n)$ we thus find that
\begin{equation*}
\begin{split}
M_2 \left(\kappa; \sqrt{\frac{v}{2}} u \right) & =  M_2 \left(\kappa; \sqrt{\frac{v}{2}} \left(2 \sqrt{a_1} n_1 + \frac{a_2}{\sqrt{a_1}} n_2\right), m n_2 \right) \\
& = - \frac{1}{\pi^2} e^{-2 \pi v Q(n)} \int_{\R^2} \frac{e^{ \frac{-\pi v \omega_1^2 - \pi v \omega_2^2}{2}}}{(\omega_2 - im n_2) (\omega_1 - \kappa \omega_2 - 2 \sqrt{a_1} i n_1)} d\omega_1 d\omega_2.
\end{split} 
\end{equation*}
Letting $\omega_1 \rightarrow 2 \sqrt{a_1} \omega_1 + \kappa \omega_2$ yields the integral as
\begin{equation*}
\begin{split}
 - \frac{1}{\pi^2} e^{-2 \pi v Q(n)} \int_{\R^2} \frac{e^{ \frac{-\pi v (2 \sqrt{a_1} \omega_1 + \kappa \omega_2)^2 - \pi v \omega_2^2}{2}}}{(\omega_2 - im n_2) (\omega_1 - i n_1)} d\omega_1 d\omega_2.
\end{split} 
\end{equation*}
Then shifting $\omega_2 \rightarrow m \omega_2$ gives
\begin{equation*}
\begin{split}
& - \frac{1}{\pi^2} e^{-2 \pi v Q(n)} \int_{\R^2} \frac{e^{ \frac{-\pi v (2 \sqrt{a_1} \omega_1 + \kappa m \omega_2)^2 - \pi v m^2\omega_2^2}{2}}}{(\omega_2 - i n_2) (\omega_1 - i n_1)} d\omega_1 d\omega_2 \\
= &  - \frac{1}{\pi^2} e^{-2 \pi v Q(n)} \int_{\R^2} \frac{e^{ - 2 \pi v Q(\omega)}}{(\omega_2 - i n_2) (\omega_1 - i n_1)} d\omega_1 d\omega_2.
\end{split} 
\end{equation*}
Therefore we have that
\begin{equation*}
\begin{split}
\lim\limits_{r \rightarrow \infty} \sum_{\substack{n \in \alpha + \Z^2 \\ |n_j - \alpha_j | \leq r}} M_2 \left(\kappa; \sqrt{\frac{v}{2}} u \right)  e^{2 \pi v Q(n)}  = 
 \lim\limits_{r \rightarrow \infty} \sum_{\substack{n \in \alpha + \Z^2 \\ |n_j - \alpha_j | \leq r}} - \frac{1}{\pi^2} \int_{\R^2} \frac{e^{ - 2 \pi v Q(\omega)} }{(\omega_2 - i n_2) (\omega_1 - i n_1)} d\omega_1 d\omega_2.
\end{split}
\end{equation*}
In exactly the same fashion as \cite{HigherDepthQMFs2} we use that 
\begin{equation*}
\pi \cot(\pi x) = \lim_{r \rightarrow \infty} \sum_{k = -r}^{r} \frac{1}{x+k}
\end{equation*}
to rewrite 
\begin{equation*}
- \lim_{r \rightarrow \infty} \sum_{\substack{n \in \Z^2 \\ |n_j | \leq r}} \frac{1}{(i \omega_1 + \alpha_1 + n_1) (i \omega_2 + \alpha_2 + n_2)} = - \pi^2 \cot(\pi(i \omega_1 + \alpha_1)) \cot(\pi(i \omega_2 + \alpha_2)).
\end{equation*}
We therefore have (using Lebesgue's theorem of dominated convergence) that
\begin{equation*}
\begin{split}
\lim\limits_{r \rightarrow \infty} \sum_{\substack{n \in \alpha + \Z^2 \\ |n_j - \alpha_j | \leq r}} M_2 \left(\kappa; \sqrt{\frac{v}{2}} u \right)  e^{2 \pi v Q(n)}  = 
 \int_{\R^2} \cot(\pi i \omega_1 + \pi \alpha_1) \cot(\pi i \omega_2 + \pi \alpha_2) e^{- 2 \pi v Q(\omega) }.
\end{split}
\end{equation*}
We may then use simple trigonometric rules to split the cotangent functions into sine and cosine (and their hyperbolic counterpart) functions by use of the formula
\begin{equation*}
\cot(x+iy) = - \frac{\sin(2x)}{\cos(2x) - \cosh(2y)} + i \frac{\sinh(2y)}{\cos(2x) - \cosh(2y)}.
\end{equation*}
This gives the integral as
\begin{equation*}
\int_{\R^2} \left( \mathcal{G}_{\alpha_1} (\omega_1) \mathcal{G}_{\alpha_2} - \mathcal{F}_{\alpha_1} (\omega_1) \mathcal{F}_{\alpha_2} (\omega_2) \right) e^{- 2 \pi v Q(\omega) } d\omega_1 d\omega_2.
\end{equation*}
Next, as in \cite{HigherDepthQMFs2}, we turn to the situation when $\alpha_1 \in \Z$ and $\alpha_2 \not\in \Z$. Then there is a term in the summation where $u_1 - \kappa u_2 = 0$. However, in view of \eqref{Equation: relation between M_2 and E_2} it still makes sense to consider our function towards this locus of discontinuity of $M_2$.
We are free to assume $\alpha_1 = 0$, since it is clear that the Mordell integral is invariant under $\alpha \rightarrow \alpha_1 + 1$. Then we consider the integral
\begin{equation}\label{Equation: Mordell integral at alpha1 = 0}
\begin{split}
- 2 \lim\limits_{\alpha_1 \rightarrow 0} \int_{\R^2} \mathcal{F}_{\alpha_1} (\omega_1) \mathcal{F}_{\alpha_2} (\omega_2)  & e^{2 \pi i \tau  Q(\omega) }  d\omega_1 d\omega_2 \\
 & = - \int_{\R^2} \mathcal{F}_{0} (\omega_1) \mathcal{F}_{\alpha_2} (\omega_2)  e^{2 \pi i \tau \left( a_1 \omega_1^2  + a_3 \omega_2^2 \right) } \sum_{\pm} \pm e^{\pm 2 \pi i \tau a_2 \omega_1 \omega_2}   d\omega_1 d\omega_2,
\end{split}
\end{equation}
where by $\sum_{\pm}$ we mean the sum over possible choices of $+$ and $-$. We see that 
\begin{equation*}
\mathcal{F}_{0} (\omega_1) =  \frac{\sinh(2 \pi \omega_1)}{\cosh(2 \pi \omega_1) - 1}
\end{equation*}
has a pole at $\omega_1 = 0$. Therefore, we write 
\begin{equation*}
\mathcal{F}_0 (\omega_1) = \left(\mathcal{F}_0 (\omega_1) - \frac{1}{\pi \omega_1} \right) + \frac{1}{\pi \omega_1}.
\end{equation*}
The contribution of the first term of the left-hand side to \eqref{Equation: Mordell integral at alpha1 = 0} is then seen to be
\begin{equation*}
\begin{split}
- \int_{\R^2} \left(\mathcal{F}_0 (\omega_1) - \frac{1}{\pi \omega_1} \right) \mathcal{F}_{\alpha_2} (\omega_2) & e^{2 \pi i \tau \left( a_1 \omega_1^2  + a_3 \omega_2^2 \right) } \sum_{\pm} \pm e^{\pm 2 \pi i \tau a_2 \omega_1 \omega_2}   d\omega_1 d\omega_2 \\
 = & - \int_{\R^2} \left(\mathcal{F}_0 (\omega_1) - \frac{1}{\pi \omega_1} \right) \mathcal{F}_{\alpha_2} (\omega_2) e^{2 \pi i \tau \left( a_1 \omega_1^2  + a_3 \omega_2^2 \right) } e^{ 2 \pi i \tau a_2 \omega_1 \omega_2}   d\omega_1 d\omega_2 \\
& -  \int_{\R^2} \left(\mathcal{F}_0 (\omega_1) - \frac{1}{\pi \omega_1} \right) \mathcal{F}_{\alpha_2} (\omega_2) e^{2 \pi i \tau \left( a_1 \omega_1^2  + a_3 \omega_2^2 \right) } e^{ - 2 \pi i \tau a_2 \omega_1 \omega_2}   d\omega_1 d\omega_2.
\end{split}
\end{equation*}
Changing $\omega_1 \rightarrow - \omega_1$ in the second integral gives overall
\begin{equation*}
 = - 2 \int_{\R^2} \left(\mathcal{F}_0 (\omega_1) - \frac{1}{\pi \omega_1} \right) \mathcal{F}_{\alpha_2} (\omega_2)  e^{2 \pi i \tau Q(\omega)  }    d\omega_1 d\omega_2.
\end{equation*}
We are left to investigate the contribution arising from $\frac{1}{\pi \omega_1}$ to \eqref{Equation: Mordell integral at alpha1 = 0}. For this, we write
\begin{equation}\label{Equation: splitting F_2}
 \mathcal{F}_{\alpha_2} (\omega_2) = \left( \mathcal{F}_{\alpha_2} (\omega_2) -  \mathcal{F}_{\alpha_2} \left( \left( \omega_2 \pm \frac{a_2}{2a_3} \omega_1  \right) \right) \right) +  \mathcal{F}_{\alpha_2} \left( \left( \omega_2 \pm \frac{a_2}{2a_3} \omega_1  \right) \right).
\end{equation}
Note in particular that we introduce the arguments in the $\mathcal{F}_{\alpha_2}$ functions coming from the diagonalisation of the quadratic form 
\begin{equation*}
a_3 \omega_2^2 \pm a_2 \omega_1 \omega_2 + a_1 \omega_1^2 = a_3 \left( \omega_2 \pm \frac{a_2}{2a_3} \omega_1 \right)^2 + \left( a_1 - \frac{a_2^2}{4a_3} \right) \omega_1^2.
\end{equation*}
The first term of \eqref{Equation: splitting F_2} yields the contribution
\begin{equation*}
- \frac{2}{\pi} \int_{\R^2} \frac{1}{\omega_1} \left( \mathcal{F}_{\alpha_2} (\omega_2) -  \mathcal{F}_{\alpha_2} \left( \left( \omega_2 \pm \frac{a_2}{2a_3} \omega_1  \right) \right) \right)  e^{2 \pi i \tau Q(\omega)}   d\omega_1 d\omega_2.
\end{equation*}
The contribution of the final term is seen to be
\begin{equation*}
- \int_{\R} \frac{e^{2 \pi i \tau \left( a_1 - \frac{a_2^2}{4a_3} \right) \omega_1^2}}{\omega_1}  \int_{\R} \sum_{\pm} \pm  \mathcal{F}_{\alpha_2} \left( \left( \omega_2 \pm \frac{a_2}{2a_3} \omega_1 \right) \right) e^{2 \pi i \tau a_3 \left( \omega_2 \pm \frac{a_2}{2a_3} \omega_1 \right)^2}   d\omega_1 d\omega_2.
\end{equation*}
Inspecting the inner integral, the term with a minus sign under the change of variables $\omega_2 \rightarrow \omega_2 + \frac{a_2}{a_3} \omega_1$ is seen to cancel with the term with positive sign, thus giving overall no contribution.

The argument when $\alpha_1 \not\in \Z$ and $\alpha_2 \in \Z$ runs in a similar way, and this completes the proof.
\end{proof}

\end{document}